\documentclass[11 pt]{amsart}
\usepackage{amssymb,latexsym,amsmath,enumerate,setspace, mathtools}
\usepackage{mathrsfs}
\input{xypic}
\xyoption{all}
\numberwithin{equation}{section}

\newcommand{\R}{{\mathbb R}}
\newcommand{\Z}{{\mathbb Z}}
\newcommand{\N}{{\mathbb N}}
\newcommand{\C}{{\mathbb C}}

\newcommand{\g}{\mathfrak{g}}
\newcommand{\ki}{\mathfrak{k}}
\newcommand{\p}{\mathfrak{p}}
\newcommand{\ai}{\mathfrak{a}}

\newcommand{\szeta}{\ensuremath{Z_{\Gamma}(s, \chi)}}

\newcommand{\SO}{{\rm SO}}

\def\a{\mathfrak{a}}

\def\m{\mathfrak{m}}
\def\p{\mathfrak{p}}

\numberwithin{equation}{section}
\newtheorem{theorem}[equation]{Theorem}

\newtheorem{lemma}[equation]{Lemma}

\newtheorem{definition}[equation]{Definition}
\newtheorem{remark}[equation]{Remark}

\title[length-holonomy spectrum]{On length-holonomy spectrum of three  dimensional compact hyperbolic manifolds}
\author{Chandrasheel Bhagwat and Ayesha Fatima}

\address{Indian Institute of Science Education and Research\\ Pune\\ India.}
\email{cbhagwat@iiserpune.ac.in, ayesha.fatima@students.iiserpune.ac.in}
\date{\today}

\begin{document}
\maketitle

\begin{abstract}
In this paper we establish a strong multiplicity one type property for the length-holonomy spectrum for the three dimensional compact hyperbolic spaces. We use the analytic properties of Selberg-Gangolli-Wakayama zeta functions associated to compact hyperbolic spaces.
\end{abstract}


\section{Introduction} 

The primitive closed geodesics on a compact hyperbolic surface are analogous to primes. Selberg constructed a generalisation of the Riemann zeta function, 
called the Selberg zeta function, for compact hyperbolic surface of genus bigger than 2 in \cite{SEL} which brought the analogy of primes and geodesics to the fore.
Such a space is of the form $\Gamma \backslash \mathscr{H}$, 
where $\mathscr{H} = SL(2, \R) / SO(2)$ is the upper half plane and $\Gamma$ is a discrete subgroup of $SL(2, \R)$.  
The Selberg zeta function is a complex valued function associated to the data $(\Gamma , \, \chi)$, 
where $\chi$ is the character of a finite dimensional unitary representation $T$ of $\Gamma$.

Further generalisations of the Selberg zeta function were given by Wakayama \cite{WAK} and Gangolli \cite{GAN2}.  
Gangolli constructed the zeta function $\szeta$ for general compact locally symmetric space of negative curvature 
$X_{\Gamma}$, 
and showed how the location and the order of the zeros of $\szeta$ gives information about the
spectrum of the Laplace-Beltrami operator of $X_{\Gamma}$ and about the topology of $X_{\Gamma}$.  

In this paper we establish a strong multiplicity one type property for the length-holonomy spectrum for the three dimensional compact hyperbolic spaces. The strong multiplicity one theorem is a classical theorem in modern number theory (see \cite{LANG}). Some of its analogues in the context of Lie groups and associated  locally symmetric spaces are established  in \cite{BR}, \cite{BR2}, \cite{KEL}.

We briefly describe the organisation of rest of the paper. 
In  Section \ref{lhs}, we give the relevant definitions needed to describe 
spaces of the type $X_{\Gamma}$ and review some basic calculations in the context of the group $\SO(3,1)^{\circ}$   and its Lie algebra ${\mathfrak so}(3,1)$. In Section \ref{zeta}, we review the theory of various Zeta functions associated to compact locally symmetric spaces. In last Section \ref{main-result}, we describe our main result and its proof.

\bigskip

\section{Length-Holonomy spectrum for compact hyperbolic spaces}\label{lhs}

\subsection{Symmetric and locally symmetric spaces}
Let $M$ be Riemannian manifold.  
Let $p$ be a point in $M$ and $N_0$ be a symmetric neighbourhood of $0$ in $T_pM$ (the tangent space of $p$). Let $N_p$ be $\exp_{p}(N_0)$, where $\exp_p$ is the exponential map from $T_pM$ to $M$.
For any $q \, \in \, N_p$, we consider the geodesic 
$t \, \longrightarrow \, \gamma(t)$ within $N_{p}$ such that $\gamma(0) \, = \, p$ and $\gamma(1) \, = \, q$.
The mapping $q \, \longmapsto \, \gamma(-1)$ of $N_p$ onto itself is called 
\text{geodesic symmetry} with respect to the point $p$.

\begin{definition}[Riemannian Locally Symmetric Space]
A Riemannian manifold $M$ is called a Riemannian locally symmetric space if 
for each $p\, \in \,  M$ there exists a normal neighbourhood 
of  $p$ on which the geodesic symmetry with respect to $p$ is an isometry.  
\end{definition}

\begin{definition}[Riemannian Globally Symmetric Space]
An analytic Riemannian manifold is called globally symmetric if 
each point $p \, \in \, M$ is the fixed point of an involutive isometry $s_p$ on $M$.   
\end{definition}

It is known that
(\cite[Ch. IV, Lemma 3.1]{HEL}) for each $p\, \in \, M$ 
there exists a normal neighbourhood $N_p$ of $p \in \,M$ such that $s_p$ is 
the geodesic symmetry on $N_p$.  
It is also known that locally symmetric Riemannian manifolds, which are not globally symmetric, 
can be constructed as quotients of Riemannian globally symmetric spaces by 
discrete groups of isometries with no fixed points.

For any Riemannian manifold $M$, let $I(M)$ be the group of isometries of $M$, endowed with the compact open topology.
The identity component of $I(M)$ is denoted by $I(M)^{\circ}$.  
It is known that, when $M$ is a globally symmetric space, $I(M)$ has a Lie group structure.  
The following theorem (\cite[Ch. IV, Theorem 3.3]{HEL}) gives the group theoretic description of 
locally and globally symmetric Riemannian spaces.  

\begin{theorem}
Let  $M$ be a Riemannian globally symmetric space and $p_0$ be any point in $M$.  
If $G \, = \, I(M)^{\circ}$, and $K$ is the subgroup of $G$ which leaves $p_0$ fixed, 
then $K$ is a compact subgroup of the connected group $G$ and 
$G/K$ is analytically diffeomorphic to $M$.    
\end{theorem}

Using the above theorem, one can represent any compact locally symmetric space of negative curvature $X_{\Gamma}$ as 
$ \Gamma \backslash G / K$, where  $G$ is a connected semisimple Lie group with finite centre, 
$K$ is a maximal compact subgroup and 
$\Gamma$ is a torsion-free uniform lattice in $G$, i.e., a co-compact torsion free discrete subgroup of $G$.  

Let $\g$ and $\ki$ be the Lie algebras of $G$ and $K$ respectively.  
Let $\g \, = \, \ki \, \oplus \, \p$ be the Cartan decomposition of $\g$ 
with respect to the involution $\theta$ determined by $\ki$.  
More specifically, $\ki$ and $\p$ are the $+1$ and $-1$ eigen-spaces of the involution map $\theta$.  
Let $\ai_{\p}$ be the maximal abelian subalgebra of $\p$.  
The subalgebra $\ai_{\p}$ can be extended to an algebra $\ai$ 
such that $\ai \, = \, \ai_{\p} \, + \, \ai_{\ki}$, where 
$\ai_{\ki} \, = \, \ai \, \cap \, \ki$ and $\ai_{\p} \, = \, \ai \, \cap \, \p$.  
The \textit{real rank} of $G$ is defined to be the dimension of $\ai_{\p}$, $\text{dim}(\ai_{\p})$.
The rank of a locally symmetric Riemannian manifold is defined to be the rank of $G$.
It is assumed all the locally symmetric Riemannian manifolds under consideration are of real rank $1$, i.e., 
$\text{dim}(\ai_{\p}) \, = \, 1$.

Let $\g^{\C}$ and $\ai^{\C}$ denote the complexifications of $\g$ and $\ai$ and 
let $\Phi(\g^{\C}, \ai^{\C})$ denote the set of roots of $(\g^{\C}, \ai^{\C})$.
An element of $\alpha$ in the dual space of  $\ai^{\C}$ is called a \textit{root} if 
$$L_{\alpha} \, = \lbrace g \, \in \, \g^{\C} \, |\quad  [h, g]  = \alpha(h)g \quad
\forall \, h \, \in \, \ai^{\C} \rbrace$$ 
is a non-zero subspace of $\g^{\C}$.  

\begin{definition}[Ordered Vector Space]
Let $V$ be a finite-dimensional vector space  over $R$.
$V$ is said to be an \textit{ordered vector space} if it is an ordered set and 
the ordering relation $>$ satisfies  the  conditions: 
\begin{itemize}
 \item $X \, > \, 0$, $Y \, > \, 0$ implies that $X \, + \, Y \, > 0$.
 \item  If  $X \, > \, 0$ and $a$ is a positive  real number, then $aX \, > \, 0$.
\end{itemize}
\end{definition}
If $\lbrace \, X_1, \, \dots \, , \, X_n \rbrace$ is a basis of $V$, 
then $V$ can be turned into an ordered vector space by saying $X \, > \, 0$ 
if $X \, = \, \sum\limits_{i = 1}^{n} a_i X_i$ 
and the first non-zero number in the sequence $a_1, \, \ldots , \, a_n$ is $> \, 0$.

Suppose $W$ a subspace of $V$.  
Let $V^{*}$ and $W^{*}$ denote their duals.  Considered them turned into ordered vector spaces.  
The orderings are said to be \textit{compatible}, if $\lambda \in V^*$ is positive whenever its restriction 
$\bar{\lambda}$ to $W$ is positive.  

We choose compatible ordering in the dual spaces of $\ai_{\p}$ and $\ai_{\p}+i \ai_{\ki}$, respectively.
Since each root $\alpha \in \Phi$ is real valued on $\ai_{\p}+i \ai_{\ki}$, this gives an ordering of $\Phi$.  
Let $\Phi^{+}$ be the set of positive roots under this order.  Let
$$P_{+} = \lbrace \alpha \in \Phi^{+} \mid \alpha \not\equiv 0 \, ~\text{on}~ \, \ai_{\p}\rbrace,$$
$$P_{-} = \lbrace \alpha \in \Phi^{+} \mid \alpha \equiv 0 \, ~\text{on}~ \, \ai_{\p}\rbrace.$$

Let $\Sigma$ be the set of restrictions of elements of $P_{+}$ to $\a_{\p}$.  
Since $ {\rm dim} (\a_{\p}) = 1$, it is known that 
we can choose $\beta \in \Sigma$ such that $2\beta$ is the only other possible element in $\Sigma$.
Let $p$ be the number of roots in $P_{+}$ whose restriction to $\ai_{\p}$ is $\beta$ and let 
$q$ be the number of roots in $P_{+}$ whose restriction is $2\beta$.
We choose $H_{0} \in a_{\p}$ such that $\beta(H_{0}) =1$.  Let $\rho = \dfrac{1}{2} \sum \limits_{\alpha \in P_{+}} \alpha$ be the half-sum of roots in $P_{+}$.  
We denote the number $\rho(H_{0})$ by $\rho_{0}$.  \smallskip

If $X_{\alpha}$ is the root vector belonging to $\alpha$, for any $\alpha \, \in \, \Phi^{+}$, then we put 
$\mathfrak{n}^{\C} \, = \, \bigoplus \limits_{\alpha \, \in \, P^{+}} \C X_{\alpha}$.  
If $\mathfrak{n} \, = \, \mathfrak{n}^{\C} \, \cap \, \g$, then 
$\g \, = \, \ki \, + \, \ai_{\p} \, + \, \mathfrak{n}$ is the \textit{Iwasawa decomposition} of $\g$.  
Furthermore, if we let $A_{\p} \, = \, \text{exp}(\ai_{\p})$ and $N \, = \, \text{exp}(\mathfrak{n})$, 
then we get $G \, = \, K A_{\p} N$, which is the 
\textit{Iwasawa decomposition of $G$}.  
Also, we can write $G \, = K A_{\p} K$, which is called the \textit{Cartan decomposition} of $G$.
The Weyl group of $(G, \, A_{\p})$ is denoted by $\mathcal{W}$.
The real dual of $\ai_{\p}$ is denoted by $\Lambda$, while the complexification 
is denoted by $\Lambda^{\C}$.  

\bigskip

\subsection{Holonomy and length of geodesics}

In this section, we describe the holonomy class of a closed geodesic in $X_{\Gamma} \, = \, \Gamma \backslash G/K$ 
in the case where $G \, = \, \SO(n, \, 1)^{\circ}$, the connected component of identity in 
the isometry group $\SO(n,\,  1)$ of the hyperbolic space $\mathbb{H}_n$.
We then give the definition of the length-holonomy spectrum of $X_{\Gamma}$.

Let $\lambda \, : \, [0, \, 1] \, \longrightarrow \, X$ be a closed geodesic in $X_{\Gamma}$.
Let $\lambda(0) \, = \, \lambda(1) \, = \, P$ and $\lambda^{'}(0) \, = \, v$.  \medskip

\begin{definition}[Parallel Transport of a vector along a geodesic]
Parallel transport of a vector $v$ along a geodesic $\lambda$ in $X_{\Gamma}$ is a 
vector field $Y_v$ along $\lambda$ in $X_{\Gamma}$ such that
\begin{enumerate}
 \item $Y_v(t) \, \in \, T_{\lambda(t)}X_{\Gamma}$, 
 \item $Y_v(t)$ is parallel to $v$ for all $t \, \in \, [0, \, 1]$, 
 \item $Y_v(0) \, = \, v$
\end{enumerate}
\end{definition} \medskip

Consider the orthogonal complement $W$ of $v =  \, \lbrace w \, \in T_{P}X_{\Gamma}:\, \langle v, \, w \rangle \, = \, 0 \rbrace $ in $T_P X_{\Gamma}$.
Using the parallel transport define the map 
\begin{align*}
T \, : W \, &\longrightarrow \, W \\
w \, &\longmapsto \, Y_w(1)
\end{align*}

We now fix an orthonormal basis $\lbrace v_1, \, v_2, \, \ldots , \, v_n-1 \rbrace$ of $W$ 
which extends to the basis $\lbrace v_1, \, v_2, \, \ldots , \, v_n-1, \, v \rbrace$ of $T_P X_{\Gamma}$ such that ${\rm det} [v_1, \, v_2, \, \ldots , \, v_n-1, \, v] = 1$.  
Then it can be checked that $T(v) \, = \, v$ and $ \langle T(w_1), \, T(w_2) \rangle  \, = \, \langle w_1, \, w_2 \rangle$ 
for all $w_1, \, w_2 \, \in \, W$.  
Therefore the map $T$, which depends only on the geodesic $\lambda$, defines an element $A$ of 
$\SO(n \, - \, 1)$ up to conjugacy.  \medskip

\begin{definition}[Holonomy Class]
The conjugacy class of $A$ in $SO(n-1)$, as defined above, is called the 
holonomy class of the closed geodesic $\lambda$ and is denoted by $h_{\lambda}$.
\end{definition}
For any $n$, let the set of conjugacy classes of $\SO(n)$ be denoted by $ \mathcal{M}_{n}$. When $G \, = \, \SO(n,1)^{\circ}$, it can be checked that $M$, the centraliser of $A_{\p}$ in $K$ is 
equal to $\SO(n-1)$.  
Furthermore, it is known that the holonomy class of geodesic $\lambda$ is 
precisely the conjugacy class $[m_{\lambda}]$.\medskip

\begin{definition}[Length-holonomy Spectrum]
The length-holonomy spectrum of $X_{\Gamma}$ is defined to be 
the function $\mathfrak{L}_{\Gamma}$ defined on $\R \times \mathcal{M}_{n \, - \, 1}$ by, 
$$\mathfrak{L}_{\Gamma}(a, \, [M]) ~ = ~\text{ number of conjugacy classes}~ [\gamma] ~\text{ in} \Gamma ~\text{ with }
(l(\gamma), \, h_{\gamma}) = (a, [M])$$
\end{definition}

\bigskip

\subsection{The Lie group $\SO(3, 1)^{\circ}$ and Lie algebra $\mathfrak{so}(3, 1)$ }\label{so-3-1}

We now consider the special case of three dimensional hyperbolic space . Here the connected component of isometry group is $G = \SO(3, 1)^{\circ}$. First we write an explicit description of the Iwasawa decomposition.  

The Lie algebra of $\SO(3,1)^{\circ}$ denoted by $\mathfrak{so}(3,1)$ is given as 
$$\{ A \, \in \, M_{4}(\R) \,  \mid \, A^{\intercal} J \, + \, JA = 0 \}$$ where 
$$
  J=
  \left( {\begin{array}{cccc}
   1 & 0 & 0 & 0 \\
   0 & 1 & 0 & 0 \\
   0 & 0 & 1 & 0 \\
   0 & 0 & 0 & -1 \\
  \end{array} } \right)
$$
It can be checked that $$\mathfrak{so}(3, 1) = \left\lbrace 
\left( {\begin{array}{cc}
   B & u \\
   u^{\intercal} & 0 
  \end{array} } \right) \, \middle| \, u \, \in \R^3 \, \text{and}\, B^{\intercal} = -B, \, B \, \in \, M_{3}(\R) \right\rbrace $$
Thus any matrix of $\mathfrak{so}(3, 1)$ can be uniquely written as
$$ \left( {\begin{array}{cc}
   B & u \\
   u^{\intercal} & 0 
  \end{array} } \right) =
  \left( {\begin{array}{cc}
   B & 0 \\
   0 & 0 
  \end{array} } \right) +
  \left( {\begin{array}{cc}
   0 & u \\
   u^{\intercal} & 0 
  \end{array} } \right)$$
where the first matrix is skew symmetric and the second matrix is symmetric 
and both belong to $\mathfrak{so}(3, 1)$.

Define $\ki$ and $\p$ as follows: 
$$\mathfrak{k} = \left\lbrace 
\left( {\begin{array}{cc}
   B & 0 \\
   0 & 0 
  \end{array} } \right) \, \middle| \, B^{\intercal} = -B \right\rbrace$$  
$$\mathfrak{p} = \left\lbrace 
\left( {\begin{array}{cc}
   0 & u \\
   u^{\intercal} & 0 
  \end{array} } \right) \, \middle| \, u \, \in \, \R^3 \right\rbrace.$$\\
  
Both $\ki$ and $\p$ are subspaces of $\g = \mathfrak{so}(3,1)$ and $\ki$ is a sub-algebra of $\g$.  
In fact, $\ki$ is isomorphic to $\mathfrak{so}(3)$.
Furthermore, we have $[\ki, \ki] \subset \ki, \, [\p, \p] \subset \ki$ and $[\ki, \p] \subset \p$.  
We have the decomposition $\mathfrak{so}(3, 1) \, = \ki \, \oplus \, \p$, which gives the 
Cartan Decomposition of $\mathfrak{so}(3, 1)$.  
Also, if $d\Theta$ is the derivative of Cartan involution map $\Theta(A) =\, (A^{\intercal})^{-1}$, 
then $\ki$ and $\p$ are the $+1$ and the $-1$ eigenspaces of $d\Theta$ respectively.

Let $\ai_{\p}$ be the matrices of the form 
$$
  B=
  \left( {\begin{array}{cccc}
   0 & 0 & 0 & 0 \\
   0 & 0 & 0 & 0 \\
   0 & 0 & 0 & \alpha \\
   0 & 0 & \alpha & 0 \\
  \end{array} } \right)
$$

Then $\ai_{\p}$ is an abelian subalgebra of $\p$.  
It can be checked that it is in fact the maximal abelian subalgebra in $\p$.  

Consider the matrices in $\mathfrak{so}(3,1)$ of the form 
$$\left( {\begin{array}{cccc}
   0 & 0 & -a & a \\
   0 & 0 & -b & b \\
   a & b & 0 & 0 \\
   a & b & 0 & 0 \\
  \end{array} } \right)
$$
where $a, \, b \, \in \, \R$.  
Such matrices form an abelian subalgebra $\mathfrak{n}$ of $\mathfrak{so}(3,1)$.
Furthermore, $\mathfrak{so}(3,1) \, = \, \ki \, \oplus \, \ai_{\p} \, \oplus \, \mathfrak{n}$.  
This gives the Iwasawa Decomposition of $\mathfrak{so}(3,1)$.  

Define $\mathfrak{m} = \left\lbrace 
\left( {\begin{array}{cc}
   B & 0 \\
   0 & 0 
  \end{array} } \right) \, \in \, M_{4}(\R) \middle| \, B \, \in \, \mathfrak{so}(2) \right\rbrace$.
It can be checked that $\mathfrak{m}$ is a subalgebra of $\ki$ and 
is the centraliser of $\ai_{\p}$ in $\mathfrak{so}(3,1)$.  
Let $\ai_{\ki} $ be defined by

$$
\ai_{\ki} = \left\lbrace 
\left( {\begin{array}{cccc}
   0 & b & 0 & 0 \\
  -b & 0 & 0 & 0 \\
   0 & 0 & 0 & 0 \\
   0 & 0 & 0 & 0
  \end{array} } \right) \, \middle| \, b \, \in \, \R \right\rbrace.
$$

Then $\ai_{\ki}$ is a maximal abelian subalgebra of $\ki$.  
Also, $\ai_{\ki}$ commutes with $\ai_{\p}$.
Therefore $\ai \, = \, \ai_{\ki} \, + \, \ai_{\p}$ is a maximal abelian subalgebra of $\g$ 
such that $\ai \cap \p = \ai_{\p}$ and $\ai \cap \ki = \ai_{\ki}$.
Then the following $\ai $ defines a Cartan subalgebra of $\mathfrak{so}(3,1)$.
$$\ai = \left\lbrace 
\left( {\begin{array}{cccc}
   0 & b & 0 & 0 \\
  -b & 0 & 0 & 0 \\
   0 & 0 & 0 & \alpha \\
   0 & 0 & \alpha & 0
  \end{array} } \right) \, \middle| \, \alpha, \, b \, \in \, \R \right\rbrace
$$

We can also check that 
$$
\exp{\ai} \, = \, A \, = \left\lbrace 
\left( {\begin{array}{cccc}
   \cos{b} & \sin{b} & 0 & 0 \\
  -\sin{b} & \cos{b} & 0 & 0 \\
   0 & 0 & \cosh{\alpha} & \sinh{\alpha} \\
   0 & 0 & \sinh{\alpha} & \cosh{\alpha}
  \end{array} } \right) \,  \right\rbrace,
$$
$$
\exp{\ai_{\p}} \, = \, A_{\p} \, = \left\lbrace 
\left( {\begin{array}{cccc}
   1 & 0 & 0 & 0 \\
   0 & 1 & 0 & 0 \\
   0 & 0 & \cosh{\alpha} & \sinh{\alpha} \\
   0 & 0 & \sinh{\alpha} & \cosh{\alpha}
  \end{array} } \right) \,  \right\rbrace,
$$
$$
\exp{\m}\, = \, M \, = \left\lbrace 
\left( {\begin{array}{cccc}
   \cos{b} & \sin{b} & 0 & 0 \\
  -\sin{b} & \cos{b} & 0 & 0 \\
   0 & 0 & 0 & 0 \\
   0 & 0 & 0 & 0
  \end{array}} \right) \,  \right\rbrace.
 $$

Therefore, in this case, $M \, \cong \, \SO(2)$ and 
hence the holonomy class of a closed geodesic $\lambda$ is given by a parameter $b \, \in \, (0, \, 2\pi]$.  
The parameter $\alpha$ gives the length of the closed geodesic.

Let $\g$ and $\ai$ be as above.  
Let $\g^{\C}$ and $\ai^{\C}$ denote the complexifications of $\g$ and $\ai$ and 
let $\Phi(\g^{\C}, \ai^{\C})$ denote the set of roots of $(\g^{\C}, \ai^{\C})$.

We now explicitly compute the roots of $\mathfrak{so}(3, 1, \, \C)$.  
(cf. Chapter IV, section 1 of \cite{KNP}.)\\

An element of $\alpha$ in the dual space of  $\ai^{\C}$ is called a root if 
$$L_{\alpha} \, = \lbrace g \, \in \, \g^{\C} \, |~ [h,g] = hg \, - \, gh = \alpha(h)g \quad
\forall \, h \, \in \, \ai^{\C} \rbrace$$ 
is a non-zero subspace of $\g^{\C}$.  
Any $h \in \ai^{\C}$ can be written as 
$$
\left( {\begin{array}{cccc}
   0 & ih_{1} & 0 & 0 \\
  -ih_{1} & 0 & 0 & 0 \\
   0 & 0 & 0 & ih_{2} \\
   0 & 0 &  ih_{2}  & 0
\end{array} } \right)
$$
for some complex numbers $h_{1}$ and $h_{2}$.
We know that $L_{\alpha} \, = \, \C E_{\alpha}$ where $E_{\alpha} \, \in \, \g^{\C}$ is of the form
$$
\left( {\begin{array}{cccc}
   0 & 0 & w & x \\
  0 & 0 & y & z \\
   -w & -y & 0 & ih_{2} \\
   x & z & ih_{2} & 0
\end{array} } \right)
$$
for some $w, \, x, \, y, \, z \, \in \, \C$.  
Therefore, $\alpha$ is a root if there exists a non-zero solution for $(w, \, x, \, y, \, z)$ 
such that $hE_{\alpha} \, - \, E_{\alpha}h \, = \, \alpha(h)E_{\alpha}$.\\

It can be shown that such a solution exists only when 
$$\alpha(h) \, \in \, \lbrace h_{1} \, + \, ih_{2}, \, h_{1} \, - \,  ih_{2}, \, -h_{1} \, - \,  ih_{2}, \, -h_{1} \, + \,  ih_{2}\rbrace$$
Denote these four roots by $\lbrace \alpha_1, \, \alpha_2, \, \alpha_3, \, \alpha_4\rbrace$, respectively.

We choose compatible ordering in the dual spaces of $\ai_{\p}$ and $\ai_{\p}+i \ai_{\ki}$, respectively.
Any element of $\ai_{\p}$ is of the form 
$$
x = \left( {\begin{array}{cccc}
   0 & 0 & 0 & 0 \\
   0 & 0 & 0 & 0 \\
   0 & 0 & 0 & a \\
   0 & 0 & a & 0
  \end{array}} \right)
$$
where $a \, \in \, \R$.
The vector space $\ai_{\p}$ and hence it's dual $\ai_{\p}^*$ are one dimensional over $\R$.  
Let $e_{1} \in \ai_{\p}^*$ such that $e_{1}(x) \, = \, a$ be a basis of $\ai_{\p}^*$.   
Any element of $\ai_{\p} + i \ai_{\ki}$ is of the form 
$$ y \, = \, \left( {\begin{array}{cccc}
   0 & ib & 0 & 0 \\
  -ib & 0 & 0 & 0 \\
   0 & 0 & 0 & a \\
   0 & 0 & a & 0
  \end{array} } \right)
$$
where $a, \, b \, \in \, \R$.  
We can extend the basis $e_1$ of $\ai_{\p}^*$ to a basis $\lbrace \, e_1, \, e_2\rbrace$ 
of ${\ai_{\p} + i\ai_{\ki}}^*$ by letting $e_2(y) \, = b$.  
This choice of basis gives a compatible ordering of the dual spaces of $\ai_{\p}$ and $\ai_{\p} + i\ai_{\ki}$.

For $y \, \in \, \ai_{\p} \, + \, i\ai_{\ki}$ as above,

$$ \alpha_1(y) \, = b \, + \, a = \, (e_1 \, + \, e_2)(y), $$
$$ \alpha_2(y) \, = b \, - \, a = \, (-e_1 \, + \, e_2)(y), $$ 
$$ \alpha_3(y) \, = -b \, - \, a = \, (-e_1 \, - \, e_2)(y), $$ 
$$ \alpha_4(y) \, = -b \, + \, a\ = \, (e_1 \, - \, e_2)(y).$$

Therefore $\alpha_1$ and $\alpha_4$ are positive roots according to the ordering fixed above.   Also, we know that for $G \, = \, \SO(2n \, + \, 1, \, 1)^{\circ}$, $p \, = 2n$ and $q \, = \, 0$.  
Therefore for $G \, = \, \SO(3, \, 1)^{\circ}$, $\rho_{0} \, = \, 1$.  

\bigskip

\section{Zeta Functions for compact locally symmetric spaces}\label{zeta}

\subsection{Uniform lattices in semisimple Lie group}
Let $G$ be a connected rank one semisimple Lie group with finite centre, $K$ be a maximal compact subgroup and 
$\Gamma$ be a torsion free uniform lattice in $G$, i.e., a co-compact torsion free discrete subgroup of $G$.  
Let $\chi$ be the character of some finite dimensional unitary representation of $\Gamma$.  
In \cite{GAN1}, Gangolli defined a zeta function $\szeta$ associated to the data $(\Gamma, \, \chi)$.  
In this section, we describe this zeta function, called the Selberg-Gangolli zeta function, and its properties.

Since $\Gamma \backslash G$ is compact, $\Gamma$ contains no parabolic elements.  
Hence, any element $\gamma \in \Gamma \setminus \lbrace e \rbrace$ is semisimple 
and therefore lies in a Cartan subgroup of $G$.  Upto conjugacy, there are only two Cartan subgroups of $G$, 
$K$ and $MA_{\p}$, where $M$ is the centraliser of $A_{\p}$ in $K$.  
Here $K$ is compact and $MA_{p}$ is non-compact.  
It is also assumed that the uniform lattice $\Gamma$ is such that it has no elements of finite order.  
Since an element is elliptic if and only if it is of finite order, 
any $\gamma \, \in \, \Gamma$ is hyperbolic and is therefore conjugate to an element of $MA_{\p}$.  
Let $A_{\p}^{+}$ be the set $\lbrace \text{exp}(tH_{0}) \, ; \, t \, \geq \, 0\rbrace$.  
Furthermore, it can be shown that it can be chosen to be conjugate to an element of $MA_{\p}^{+}$.
Let $h(\gamma) \, = \, m_{\gamma}(\gamma)h_{\p}(\gamma)$ be an element of $MA_{\p}^{+}$ conjugate to $\gamma$.  

For any $h \, \in \, A_{\p}$, let $u(h) \, = \, \beta(\text{log} \, h)$.  
Then $u \, = \, u(h)$ can be considered as a parametrisation on $A_{\p}$ via which $A_{\p}$ can be identified with $\R$.
Under this parametrisation, $A_{\p}^{+}$ corresponds to the positive real axis.  
Define $u_{\gamma} \, := \, \beta(log \, h_{\p}(\gamma))$.  
It is known that the value $u_{\gamma}$ is independent of the choice of $h(\gamma)$.  
\begin{lemma}
For any negatively curved Riemannian manifold, and hence for $X_{\Gamma}$, 
there is a bijective correspondence between the set of closed geodesic classes in 
$X_{\Gamma}$ and the set of conjugacy classes of $\Gamma$.  
\end{lemma}
Let $C_{\Gamma}$ be a set of representatives of the $\Gamma$-conjugacy classes.
\begin{definition}
For $\gamma \, \in \, \Gamma$, the length $\ell(\gamma)$ of the conjugacy class $[\gamma]$ is defined 
to be the length of the unique closed geodesic in the corresponding free homotopy class in $X_{\Gamma}$. (this is well defined as the manifold is compact and has constant negative curvature).
\end{definition}
\begin{lemma}[\cite{GAN2}]
For any $\gamma \, \neq \, 1$ in $\Gamma$, we have $\ell(\gamma) \, = \, u_{\gamma}$.  
\end{lemma}

The length spectrum of $X_{\Gamma}$ is the function 
$$ L_{\Gamma} : \R \rightarrow \mathbb{N} \cup \left\{ 0 \right\}$$
which to a real
number $l$, assigns the number of $\Gamma$-conjugacy classes $[\gamma]$ such that $\ell(\gamma) \, = \, l$.

An conjugacy class $[\gamma] \, \neq \, 1$ of $\Gamma$ is called \textit{primitive} if 
$[\gamma] \, \neq \, \, [\delta^{n}]$ for any integer $n \, > \, 1$ and $\delta \, \in \, \Gamma$.  
Furthermore, it is known that (cf. \cite{GAN2}) that any conjugacy class 
$[\gamma]$ can be written as $[\delta^{n}]$, for some primitive conjugacy class $[\delta]$ and integer $n \, \geq \, 1$.
If $[\gamma] \, = \, [\delta^{n}]$, for some primitive $[\delta]$, then the number 
$n$ is denoted by $j(\gamma)$.  
It is known that the primitive conjugacy classes correspond to primitive closed geodesics, 
which are closed geodesics which cannot be obtained by going around any other closed geodesic $n \, > \, 1$ times.
Let $P_{\Gamma}$ be the set of representatives of the primitive $\Gamma$-conjugacy classes.  

The elements of $P_{+}$ are enumerated as $\alpha_{1}, \, \alpha_{2}, \, \ldots , \, \alpha_{t} $.  
Then, we define $L$ as the semi lattice in $\ai_{\p}^{\C}$ given by 
$\lbrace \, \sum_{ i \, = \, 1}^{t} m_{i} \alpha_{i} \, ; \, m_{i} \, \geq \, 0, \, m_{i} \, \in \, \Z \rbrace$.  
For any $\lambda \, \in \, L$, the number of distinct ordered $t$-tuples $(m_1, \, m_2, \, \ldots , \, m_t)$ 
such that $\lambda \, = \, \sum_{ i \, = \, 1}^{t} m_{i} \alpha_{i}$.
The character of the Cartan subgroup $A$ which corresponds to $\lambda$ is denoted by $ \xi_{\lambda}$.  
Therefore, $\xi_{\lambda}(h) \, = \, \exp \lambda(\text{log}(h))$.\bigskip

\subsection{Selberg-Gangolli Zeta function}

For $s \, \in \, \C$, the zeta function $\szeta(s)$ is defined by first defining 
its logarithmic derivative with respect to $s$, denoted by $\Psi_{\Gamma}(s, \, \chi) $, which is 
given as the following series convergent on $ \text{Re} \, s > 2\rho_{0}$ (\cite{GAN1}):
\begin{equation}
\Psi_{\Gamma}(s, \, \chi) \, = \, \kappa \, \sum\limits_{\gamma \, \in \, C_{\gamma} \, \setminus \, \lbrace 1 \rbrace}
\chi(\gamma) \, u_{\gamma} \, j(\gamma)^{-1} \, C(h(\gamma)) \, \exp(\rho_{0} \, - \, s)\, u_{\gamma}
\end{equation}
Here $C(h(\gamma))$ is a positive function and $\kappa$ is a positive integer, 
both of which depends only on the structure of $G$
If $\szeta$ is the Selberg-Gangolli zeta function, then 

\begin{equation}
\dfrac{d}{ds} \, \log \szeta \, = \, \Psi_{\Gamma}(s, \, \chi).   
\end{equation}

Furthermore, it has been shown in \cite{GAN2} that $\szeta$ can be simplified to the following Euler product: \medskip

\begin{definition}[Selberg-Gangolli Zeta Function]

\begin{equation}
\label{szeta}
\szeta \, = \, C \, \prod\limits_{\delta \, \in \, PC_{\Gamma}} \, \prod\limits_{\lambda \, \in \, L}
({\rm det} \, ( I \, - \, T(\delta)\, \xi_{\lambda}(h(\delta))^{-1} \, {\rm exp} (-su_{\delta})))^{m_{\lambda \kappa}}, 
\end{equation}
Here, $C \, \neq \, 0$ is a constant which depends only on $G$.
\end{definition}

We state below some properties of this zeta function which have been proved in \cite{GAN2}: 
\begin{theorem}
\begin{enumerate}
\item[]
 \item $\szeta$ is holomorphic in the half plane $\text{Re} \, s \, > \, 2\rho_{0}$.  
 \item $\szeta$ has a meromorphic continuation to the whole complex plane.
 \item $\szeta$ satisfies the following functional equation 
 \begin{equation}
  Z_{\Gamma}(2 \rho _{0} - \, s, \chi) = \szeta ~\exp \int _{0} ^{s - \rho_{0}} \Phi(t) dt
 \end{equation}
where $\Phi(t) = \kappa \, {\rm vol} (\Gamma \backslash G) \chi(1) c(it)^{-1}c(-it)^{-1}$.
Here, $c(t)$ is the Harish-Chandra $c$-function, which is a function of one complex variable which depends only on 
$(G, \, K)$.
\item When $G \, = \, \SO(2n+1, \, 1)^{\circ}$ i.e.,when  the Riemannian globally symmetric space is the hyperbolic space of odd dimension, $\Phi$ is a polynomial and so is $\int _{0} ^{s - \rho_{0}} \Phi(t) dt$.  
This simplifies the functional equation.  \medskip
\end{enumerate}
\end{theorem}\bigskip

\subsection{Spherical representations}

Certain irreducible representations of the isometry group are closely related to the spectral properties of differential operators on symmetric space. In this section we review this.

\begin{definition} An irreducible unitary  representation $\pi$ of $G$ on a Hilbert space $V$
is said to be spherical if there exists a non-zero vector $v \in V$
such that
\[ \pi(k) v = v\quad \forall \ k \in K.\]
\end{definition}

We denote the set of equivalence classes of the spherical representations of $G$ by $\widehat{G}_s$.
The \textit{spherical spectrum} of a uniform lattice is the collection of spherical representations of $G$ 
which occur in the decomposition of $L^{2}(\Gamma \backslash G)$.

Let $T$ be a finite dimensional unitary representation of $\Gamma$ with characteristic $\chi$.  
We denote with $U$ the unitary representation of $G$ induced by $T$.  
Then $U$ is a discrete direct sum of irreducible unitary representations occurring with finite multiplicity.  
Let $ \left\lbrace U_j, j \geq 0 \right\rbrace$ be the spherical representations among these and let $n_j (\chi)$ be their multiplicities.
The spherical spectrum of a uniform lattice is precisely the set 
$ \left\lbrace U_j, j \geq 0 \right\rbrace$ with multiplicities 
when $T$ is chosen to be the trivial representation of $\Gamma$.

We now state two theorems of Harish-Chandra \cite{H-C} which  characterise the  irreducible unitary spherical representations, 
which are also called representations of class-1.  
\medskip

\begin{definition}[Positive definite function]\label{def1}
A complex valued continuous function $\phi$ on a topological group $G$ is called positive definite if 
$$\sum_{i,j =1}^{n} 
\phi (x_{i}^{-1} x_{j}) 
\alpha_{i} \bar{\alpha_{j}} 
\geq 0$$ for all finite sets 
$x_{1}, \, x_{2}, \, \dots x_{n}$ in $G$ and any set of complex numbers $\alpha_{1}, \, \alpha_{2}, \, \dots \alpha_{n}$.
\end{definition}
\medskip

\begin{theorem}\label{thm1}[Harish-Chandra 1]
There is a bijective correspondence between the set $\Omega$ of equivalence classes $\omega$ 
of representations of class 1 
and the set $\beta$ of all positive definite spherical functions $\varphi$ on $G$ satisfying $\varphi(e) = 1$.  
\end{theorem}
\medskip

We write $G = KA_{p} N$, Iwasawa decomposition and define $H(x)$ as the unique element in $\mathfrak{a}_{p}$ 
such that $x = k (\exp H(x)) n $, where $k \in K$ and $n \in N$.

\medskip
\begin{theorem}\label{thm2}[Harish-Chandra 2]
 Functions of the type $$\varphi _{\nu}(x) = \int_K e^{(i \nu - \rho) (H(xk))} dk , \qquad \nu \in \Lambda^{\C}$$ exhaust the 
 class of spherical functions of $G$.  
 Moreover, any two such functions $\varphi_{\nu}$ and $\varphi_{\lambda}$ are identical if and only if $\nu = \lambda ^{s}$ 
 for some $s \in \mathcal{W}$ (Weyl Group).  
\end{theorem}
\medskip

Using these two theorems, we get that each of the $U_{i}$ is uniquely determined by a spherical function $\varphi_{\nu_{j}}$, 
where the $\nu_{j} \in \Lambda^{\C}$ is uniquely determined up to action of the Weyl group.  

Let $r_{j}^{+} (\chi) = \nu_{j} (H_{0})$ and $r_{j}^{-} (\chi) = - \nu_{j} (H_{0})$.  
Put $s_j^{+} = \rho_{0} + ir_j^{+}$ and $s_j^{-} = \rho_{0} + ir_j^{-}$.  
Though these quantities depend on $\chi$, where there is no danger of confusion the $\chi$ is not explicitly mentioned.  

The following theorem from \cite{GAN1} shows that the information about the 
spherical spectrum of $X_{\Gamma}$ is encoded in $\szeta$. \medskip
 
\begin{theorem}[Zeros and Poles of $\szeta$]
\label{thm3}
\begin{enumerate}
\item[]
\item The points $s_{j}^{+}$ and $s_{j}^{-}$, with $j \geq 1$, are zeroes of $\szeta$ of order $\kappa n_{j}$ respectively.  
These are called the spherical zeroes.
\item The points $\rho_0 + ir_{k}$, $k \geq 0$ are either zeroes or poles of $\szeta$ according to whether $-\chi (1)e_{k}E$ is positive or
negative.  These are called the topological zeros or poles.  
The order of the zero or pole $\rho_{0} + ir_{k}$ is $| \chi (1) e_{k} E |$.  
\item It also follows that the $\rho_{0} + ir_k$ are either all poles or all zeroes, if they exist.  
\item Also, since the elements $r_{k}$ are purely imaginary, it follows that $\rho_{0} + ir_{k}$ equals either $-k$ or $-2k$, $k \geq 0$, 
when the $r_{k}$ exist.  
\end{enumerate} 
Here, $E$ denotes the Euler-Poincare characteristic of the manifold $X_{\Gamma}$.
\end{theorem}
When the dimension of $X_{\Gamma}$ is odd, the Euler-Poincare characteristic is zero, 
and hence the zeta function has only spectral zeros and no poles.  

The following theorem from \cite{GAN2} shows the relation between the spherical spectrum and the length spectrum.  
Simply put, it shows that the spherical spectrum of the space, and hence the Laplace spectrum, 
determines the length spectrum.
This theorem has been proved using the properties of the logarithmic derivative of $\szeta$, 
showing its use as a tool to connect the geometric and the algebraic properties of the manifold $X_{\Gamma}$.  
\begin{theorem}
Let $G$ be a connected semi-simple Lie group with finite centre 
and let $\Gamma_{1}$ and $\Gamma_{2}$ be two co-compact torsion free lattices in $G$.
Assume that the spherical spectrum of the spaces 
$X_{\Gamma_i} \, = \, \Gamma_i \ G \backslash K$ for $i = 1, \, 2$ is same.
Then, $$\lbrace l_{i} \, \in \, \R \; | \; L_{\Gamma_1}(l_i) \, \neq \, 0\rbrace \, = 
\, \lbrace l_{j} \, \in \, \R \; / \; L_{\Gamma_2}(l_j) \, \neq \, 0 \rbrace.$$
In other words, the lengths that appear in the spectrum of $X_{\Gamma_{1}}$ and $X_{\Gamma_{2}}$ are the same.
\end{theorem}\bigskip

\subsection{Selberg-Gangolli-Wakayama Zeta Function}
A further generalisation of the Selberg zeta function was defined by Wakayama in \cite{WAK}.  
This zeta function is also associated to a manifold $X_{\Gamma}$ of the type 
$\Gamma \backslash G/K$, where 
$G$ is a connected rank one semisimple Lie group with finite centre, $K$ is a maximal compact subgroup 
and $\Gamma$ is a uniform lattice in $G$.
Apart from an irreducible unitary representation of the lattice $\Gamma$, 
this zeta function also takes into consideration an irreducible unitary representation of $K$.  

For any subgroup $L$ of $G$, we denote 
the set of all equivalence classes of irreducible unitary representations of $L$
by $ \widehat{L}$.  
If $\pi \in \widehat{L}$ is a finite-dimensional representation, 
we put $\chi_{\pi} \, =\,  {\rm tr}\ \pi$ and $d_{\pi} \, = \, {\rm dim}\ \pi$.  

We recall that $M$ is the centralizer of $A_{\p}$ in $K$.  
For $\tau \in \widehat{K}$, we put 
$\widehat{M}_{\tau} = \lbrace \sigma \in \widehat{M} \, \mid \, [\sigma: \tau|_{M}] \neq 0] \rbrace$.
For ease of notation, we let $\tau_{M} := \tau |_{M}$ and $\alpha_{\sigma}\,  = \, [\sigma: \tau|_{M}]$.

The Selberg-Gangolli-Wakayama zeta function $Z_{\tau}(s)$ of $X_{\Gamma} = \Gamma \backslash G / K $, 
associated with $\tau \in \widehat{K}$ is defined by the following Euler product: \medskip

\begin{definition}[Selberg-Gangolli-Wakayama Zeta Function]
\begin{equation}
\label{wak-1}
Z_{\tau}(s) \, = \, \prod\limits_{\sigma \in \widehat{M}_{\tau}} Z_{\sigma}(s)^{[\sigma:\tau|_{M}]}
\end{equation} 
where 
\begin{equation}
\label{wak-2}
Z_{\sigma}(s) = \prod\limits_{p \in P_{\Gamma}} \prod\limits_{\lambda \in L} (1 - \chi_{\sigma}(m_{p})^{-1}
\xi_{\lambda}(h(p))^{-1}e^{-s u_{p}})^{m_{\lambda}}.
\end{equation}
\end{definition}\medskip

The Selberg-Gangolli-Wakayama zeta function satisfies many properties similar to the ones satisfied by the 
Selberg-Gangolli zeta function.  
The following theorem summarises some of the results on $Z_{\tau}$ proved 
in \cite{WAK} (Theorem 7.1-7.3, p.p.  287-291).  \medskip

\begin{theorem}
\begin{enumerate}
 \item[]
 \item $Z_{\tau}$ is holomorphic in the half plane $\text{Re} (s) \, > \, 2\rho_{0}$ and 
 has analytic continuation as a meromorphic function on the whole complex plane.
 \item $Z_{\tau}$ always have some zeros, which are called spectral zeros.  
 The location and the order of the zeros contains information about the $\tau$-spectrum $\widehat{G_{\tau}}$ 
 in $L^{2}(\Gamma \backslash G)$ where 
 $$\widehat{G_{\tau}} \, = \, \lbrace \pi \, \in \, \widehat{G};~ m_{\Gamma}(\pi) \, > \, 0, \, \tau \, \in \, \pi|_{K} \rbrace.$$
 Here $m_{\Gamma}$ is the multiplicity of a unitary representation $\pi$ of $G$ in the right regular representation 
 $\pi_{\Gamma}$ of $G$ on $L^{2}(\Gamma \backslash G)$.
 \item Apart from the spectral zeros, $\Z_{\tau}$ has certain series of zeros and poles which come from the 
 Plancherel measure.  These are called the topological zeros of $Z_{\tau}$.
 \item $Z_{\tau}$ satisfies the following functional equation: 
 \begin{equation}
  \label{wak-func}
  Z_{\tau}(2\rho_0 \, - \, s) \, = \, Z_{\tau}(s) \, 
  \text{exp}\left( \int\limits_{0}^{s \, - \, \rho_0} \, \Delta_{\tau}(t) \, dt \right).  
 \end{equation}
Here $\Delta_{\tau}$ is an expression which depends on $G, \, K, \, \chi$ and $\tau$.  
In the case when dimension of $X_{\Gamma}$ is odd, $\Delta_{\tau}$ is a polynomial and hence the functional equation simplifies to $$ Z_{\tau}(2\rho_0 \, - \, s) \, = \, Z_{\tau}(s) $$
\end{enumerate}
\end{theorem}
\medskip

\section{Towards Strong Multiplicity One Property for Length-Holonomy Spectrum for three dimensional compact hyperbolic spaces}\label{main-result}
In this section we describe a strong multiplicity one property for 
the length-holonomy spectrum of $\SO(3,1)^{\circ}$.  
More precisely, we consider two uniform lattices $\Gamma_{1}$ and $\Gamma_{2}$ of $\SO(3,1)^{\circ}$.  
The spaces associated with these lattices, $X_{\Gamma_i} = \Gamma_{i} \backslash G /K$ for $i = \, 1, \, 2$, 
have associated with them a primitive length-holonomy spectrum.

The primitive length-holonomy spectrum of $X_{\Gamma}$ is defined to be 
the function $\mathfrak{P}_{\Gamma}$ defined on $\R \, \times \, [0,2\pi)$ by,
$$\mathfrak{P}_{\Gamma}(a, \, b) \, = \# \; \text{of conjugacy classes} \; [p] \, \in \, P_{\Gamma}
\; \text{such that} \;
(a(p), \, b(p)) = (a, b).$$ 

We now define a modified length-holonomy spectrum of $X_{\Gamma}$.  
For any $p \, \in \, P_{\Gamma_i}$, we define $c(p) \, := \, \dfrac{b(p)}{a(p)}$.  
The \textit{modified primitive length-holonomy spectrum} of $X_{\Gamma}$ is defined to be 
the function $\mathfrak{M}_{\Gamma}$ defined on $\R$ by,
$$\mathfrak{M}_{\Gamma}(c) \, = \, \# \, \text{of conjugacy classes} \; [p]\, \in \, P_{\Gamma} \; 
\text{such that} \; c(p) \, = \, c.$$

We assume that the primitive length-holonomy spectra of $X_{\Gamma_1}$ and $X_{\Gamma_2}$ is same for all but 
finitely many points in $\R \times [0, \, 2\pi)$, i.e., 
$\exists$ a finite subset $S \, \subset \, \R \times [0, \, 2\pi)$ such that 
$\mathfrak{P}_{\Gamma_1}(a, \, b) \, = \, \mathfrak{P}_{\Gamma_2}(a, \, b)$ $\forall \, \in \, (a, \, b) \, \notin \, S$.\\

We first simplify the zeta function for the case $\SO(3,1)^{\circ}$.  
Consider any $p \, \in \, P_{\Gamma}$, a primitive hyperbolic conjugacy class of $\Gamma$.  
Then $p$ is conjugate to an element $h(p) = m_{p}\exp(l(p)H_{0}) \, \in \, MA_{p}^{+}$.\\

Here 
\begin{align*}
\log(h(p)) \, &= \, \left( {\begin{array}{cccc}
   0 & b(p) & 0 & 0 \\
  -b(p) & 0 & 0 & 0 \\
   0 & 0 & 0 & a(p) \\
   0 & 0 & a(p) & 0
  \end{array} } \right) \, \in \, \ai , \\[3ex]
m_{p} &= \left( {\begin{array}{cccc}
   \cos b(p) & \sin b(p) & 0 & 0 \\
  -\sin b(p) & \cos b(p) & 0 & 0 \\
   0 & 0 & 1 & 0 \\
   0 & 0 & 0 & 1
  \end{array} } \right) \, \in \, M ,\\[3ex]
\exp (l(p)H_0) \, &= \, \left( {\begin{array}{cccc}
   1 & 0 & 0 & 0 \\
   0 & 1 & 0 & 0 \\
   0 & 0 & \cosh a(p) & \sinh a(p) \\
   0 & 0 & \sinh a(p) & \cosh a(p)
  \end{array} } \right) \, \in \, A_{\p}^{+}
\end{align*}  \smallskip

for some $a(p), \, b(p) \, \in \R$.  (Here, $\ell(p) = a(p)$.)\smallskip

If $\lambda \, \in \, L$ is of the form $m_{1}(e_{1}+e_{2}) \, + \, m_{2}(e_{1}-e_{2})$, 
for $m_{1}, \, m_{2} \, \in \,  \Z^{+} \, \cup \, \lbrace 0 \rbrace$, 
then $$\lambda (\log h(p)) \, = \, (m_1 + m_2)a(p) \, + \, i(m_1 - m_2)b(p).$$
Therefore, $$\xi(h(p)) = \exp \left( (m_1 + m_2)a(p) \, + \, i(m_1 - m_2)b(p) \right).$$ \smallskip

Let $\tau$ be any irreducible representation of $\SO(3)$.
We know that any irreducible representation of $\SO(3)$ is of odd dimension and hence dimension of $\tau$ 
is of the form $2m+1$ for some m.  
Since $M = \SO(2)$ is an abelian group, any irreducible representation of $\SO(2)$ is one dimensional.  
Thus the restriction of $\tau$ to $\SO(2)$ is a direct sum of $2m+1$ one-dimensional representations, 
$$\sigma_j: \SO(2) \longrightarrow \C^{\times}$$ given by $R(\theta)\, \longmapsto \, e^{ij\theta}$, 
for $j = -m, \, \ldots , \, m$.  
Here $R(\theta) = \left( {\begin{array}{cc}
   \cos \theta & \sin \theta  \\
  - \sin \theta & \cos \theta 
  \end{array} } \right) $.  
Therefore $\chi_{\sigma_m}(m_{p}) = e^{im\theta}$.  
It can also be checked that, for $\SO(3,1)^{\circ}$, $m_{\lambda} = 1$ for any $\lambda \, \in \, L$.\medskip

With these simplifications, the Selberg-Gangolli-Wakayama Zeta function for $\SO(3,1)^{\circ}$.   is as follows: 
\begin{equation}
 \label{zeta-for-3}
 Z_{\tau}(s) \, = \,  \prod\limits_{k=-m}^{m}\prod\limits_{p \in P_{\Gamma}} \prod\limits_{\lambda \in L}
(1 - \, e^{-(ikb(p) \, + \, (m_1 + m_2)a(p) \, - \, (m_1 - m_2)b(p) \, + \, sa(p))}).
\end{equation}

We consider ratio of the corresponding zeta functions, 
\begin{equation}
\dfrac{Z_{\Gamma_1}(s)}{Z_{\Gamma_2}(s)} \, = \, \dfrac{\prod\limits_{k=-m}^{m}\, \prod\limits_{p \in P_{\Gamma_1}} \, \prod\limits_{m_1, m_2 \, \in \, \N}
(1 - \, e^{-(ikb(p) \, + \, (m_1 + m_2)a(p) \, - \, (m_1 - m_2)b(p) \, + \, sa(p))})}
{\prod\limits_{k=-m}^{m}\, \prod\limits_{q \in P_{\Gamma_2}} \, \prod\limits_{l_1, l_2 \, \in \, \N}
(1 - \, e^{-(ikb(q) \, + \, (l_1 + l_2)a(q) \, - \, (l_1 - l_2)b(q) \, + \, sa(q))})}.
\end{equation}

Under our assumption, there exist finite indexing sets $S_1$ and $S_2$ such that the above fraction becomes, 
\begin{equation}
\label{ratio-zeta}
\dfrac{Z_{\Gamma_1}(s)}{Z_{\Gamma_2}(s)} \, = \, \dfrac{\prod\limits_{k=-m}^{m}\, \prod\limits_{p \in S_1} \, 
\prod\limits_{m_1, m_2 \, \in \, \N}
(1 - \, e^{-(ikb(p) \, + \, (m_1 + m_2)a(p) \, - \, (m_1 - m_2)b(p) \, + \, sa(p))})}
{\prod\limits_{k=-m}^{m}\, \prod\limits_{q \in S_2} \, \prod\limits_{l_1, l_2 \, \in \, \N}
(1 - \, e^{-(ikb(q) \, + \, (l_1 + l_2)a(q) \, - \, (l_1 - l_2)b(q) \, + \, sa(q))})}.
\end{equation}

It can be checked that both the numerator and the denominator of the above ratio is holomorphic, which gives the following lemma: 
\begin{lemma}
The ratio $\dfrac{Z_{\Gamma_1}(s)}{Z_{\Gamma_2}(s)}$ (\ref{ratio-zeta}) is meromorphic.  
\end{lemma}

Let $T(s) \, = \, \dfrac{Z_{\Gamma_1}(s)}{Z_{\Gamma_2}(s)}$.  
We have proved that $T(s)$ is a meromorphic function.  
We also know that, for 
$i \, = \, 1, \, 2$, $Z_{\Gamma_i}$ and hence $T(s)$ admits a meromorphic continuation to $\C$.  
Therefore the expressions must match $\forall \, s\, \in \, \C$, i.e., 
\begin{multline}
\dfrac{Z_{\Gamma_1}(1\, - \, s)}{Z_{\Gamma_2}(1\, - \, s)} \, = \\ \, \dfrac{\prod\limits_{k=-1}^{1}\, \prod\limits_{p \in S_1} \, \prod\limits_{m_1, m_2 \, \in \, \N}
(1 - \, e^{-(ikb(p) \, + \, (m_1 + m_2)a(p) \, - \, (m_1 - m_2)b(p) \, + \, (1\, - \, s)a(p))})}
{\prod\limits_{k=-1}^{1}\, \prod\limits_{q \in S_2} \, \prod\limits_{l_1, l_2 \, \in \, \N}
(1 - \, e^{-(ikb(q) \, + \, (l_1 + l_2)a(q) \, - \, (l_1 - l_2)b(q) \, + \, (1\, - \, s)a(q))})}.
\end{multline}

For $G \, = \, \SO(n,1)^{\circ}$, when $n$ is odd and hence for  $G \, = \, \SO(3, \, 1)^{\circ}$, it is know that $Z_{\Gamma}(s) \, = \, Z_{\Gamma}(1\, - \, s)$ \cite{GON}. This implies that
$T(s) \, = \, T(1\, - \, s)$.  
Any zero of $T(s)$ is either a zero of the numerator in the expression for $T$ or pole of the denominator.  
But we have proved that the denominator of the expression is analytic on $\C$.  
Hence the zeros of $T(s)$ are precisely the zeros of the numerator of the expression of $T(s)$.  
Let $A$ denote the set of zeros of $T(s)$.  
It can be deduced that that the elements of $A$ are precisely the complex numbers of the form 
$s \,= \, s_1 \, + \, s_2$ where
\begin{align}
    s_1 &= -m_1\, -\, m_2\\[3ex]       
    s_2 &= \dfrac{-b(p)(m_1\, - \, m_2 \, +\, k )-2n\pi}{a(p)}.
\end{align}
\vspace{5pt}
Here $m_1 , \, m_2 \, \in \, \Z^{+}$, $k \, = \, -1, \, 0, \, 1$, $p \, \in \, S_1$ and $n \, \in \, \Z$.

We know that $T(s) \, = \, T(1-s)$.   Therefore we get that $s \, \in \, A \, \Rightarrow \, 1 \, - \, s \, \in \, A$.
But this is not true.  Therefore we get that $A \, = \, \emptyset$, i.e., $T(s)$ has no zeros.  
By a similar  we can conclude that $T(s)$ has no poles.  
Therefore the zeros of the numerator in the expression for $T(s)$ cancel out with the zeros of the denominator.  

More precisely, for any fixed $m_1, \, m_2 \, \in \, \Z^{+}$, $k_1 \, \in \, \lbrace -1, \, 0, \, 1\rbrace$, $n_1 \, \in \, \Z$
and $p \, \in \, S_{1}$, 
there exist $l_1, \, l_2 \, \in \, \Z^{+}$, $k_2 \, \in \, \lbrace -1, \, 0, \, 1\rbrace$, $n_2 \, \in \, \Z$ 
and $q \, \in \, S_2$
such that 
\begin{align}
 m_1 \, + \, m_2 \, &= \, l_1 \, + \, l_2 \label{zeroeq1} \\[3ex]
 \dfrac{-b(p)(m_1\, - \, m_2 \, +\, k_1 )-2n_1\pi}{a(p)} \, &= \, \dfrac{-b(q)(m_1\, - \, m_2 \, +\, k_2 )-2n_2\pi}{a(q)} \label{zeroeq2}
\end{align}
\vspace{5pt}

Let $A_{\Gamma_1}$ be the multiset of zeros of the expression in the numerator of $T(s)$ 
and let $A_{\Gamma_2}$ be the multiset of zeros of the expression in the denominator of $T(s)$.  
In terms of the new notation, the equations \ref{zeroeq1} and \ref{zeroeq2} state that $A_{\Gamma_1} \, = \, A_{\Gamma_2}$.  

Consider now the elements of $A_{\Gamma_1}$ which lie on the line $Re(s) \, = \, 0$.  
We have $(m_1, \, m_2) \, = \, (0, \, 0)$ for these elements.  
Also, each such zero appears with the same multiplicity in $A_{\Gamma_2}$ and will have 
$(l_1, \, l_2) \, = \, (0, \, 0)$.  

Therefore we have the following equality of the multisets, 
\begin{multline}\label{eqzeroline}\bigg\{ \dfrac{-b(p)k_1 - 2n_1\pi}{a(p)} \, 
\mid \, p \, \in \, S_1, \, k_1 \, = \, -m, 1-m \ldots, \, m-1, m, \, n_1 \, \in \, \Z \bigg\} \, = \\[2ex]
\bigg\{ \dfrac{-b(q)k_2 - 2n_2\pi}{a(q)} \, 
\mid \, q \, \in \, S_2, \, k_2 \, = \, -m, 1-m \ldots, \, m-1, m, \, n_2 \, \in \, \Z \bigg\}\
\end{multline}
\vspace{5pt}

Now consider the case when $\tau$ is the trivial representation of $\SO(3,1)^{\circ}$.  For this $\tau$, we have $m \, = \, 0$ and \ref{eqzeroline} becomes the following equality of multisets:\medskip
\begin{multline}
\label{res1eq}
\bigg\{ \dfrac{2n_1\pi}{a(p)} \, \mid \, p \, \in \, S_1, \, n_1 \, \in \, \Z \bigg\} \, = 
\bigg\{ \dfrac{2n_2\pi}{a(q)} \, \mid \, q \, \in \, S_2, \, n_2 \, \in \, \Z \bigg\}\
\end{multline} \medskip

Let $s_0$ the element in the set on the left hand side which is strictly positive and closest to $0$.  
For $n \, \leq \, 0$, $\dfrac{2n\pi}{a(p)} \, \leq \, 0$.  
For $n \, > \, 1 $,  $\dfrac{2n\pi}{a(p)} > \dfrac{2\pi}{a(p)}$.  
Therefore, $s_0 \, = \, \dfrac{2\pi}{a(p)}$ for some $p \, \in \, S_1$.  
Let $s^{'}_{0}$ be the element in the set on the right hand side which is strictly positive and closest to $0$.  
By similar argument, it can be concluded that $s_{0}^{'} \, = \, \dfrac{2\pi}{a(q)}$ for some $q \, \in \, S_2$.
Now, the equality of the sets in equation \ref{res1eq} implies that for the $p \, \in \, S_1$ corresponding to $s_{0}$, 
there exists a $q \, \in \, S_2$ (corresponding to $s_{0}^{'}$) such that 
$$ \dfrac{2\pi}{a(p)} \, = \, \dfrac{2\pi}{a(q)}$$

Suppose the zero $\dfrac{2\pi}{a(p)}$ has multiplicity $ > \, 1$.  
It has the same multiplicity in the set on right hand side of eq. \ref{res1eq}.  
Therefore, there exist a $0 \, < \, n^{'} \, \in \, \Z$ and $q^{'} \, \in \, S_2$ such that 
$\dfrac{2\pi}{a(p)} \, = \, \dfrac{2n^{'}\pi}{a(q^{'})}$.  
If $n^{'} \, > \, 1$, then 
$$ \dfrac{2\pi}{a(q^{'})} \, < \, \dfrac{2n^{'} \pi}{a(q^{'})} \, = \, \dfrac{2\pi}{a(q)}.$$
This is a contradiction.  Hence $n^{'} \, = \, 1$ and  
$$ \dfrac{2\pi}{a(q^{'})} \, = \, \dfrac{2 \pi}{a(q^{'})} \, = \, \dfrac{2\pi}{a(q)}.$$
Therefore the zeros $\dfrac{2\pi}{a(p)}$ and $\dfrac{2\pi}{a(q)}$ appear with the same multiplicity 
in the sets on the left hand side and the set on the right hand side of eq. \ref{res1eq}, respectively.  

We can now remove the points of the form $\dfrac{2n\pi}{a(p)}$  
and points of the form $\dfrac{2n\pi}{a(q)}$ from the sets on the left-hand side and right-hand side 
of eq.  \ref{res1eq} respectively and repeat the argument to get the following result: \medskip

\begin{theorem}
\label{res1}
Let $G \, = \, \SO(3,1)^{\circ}$ and $\Gamma_1$ and $\Gamma_2$ be two uniform lattices in $G$ such that 
$\mathfrak{P}_{\Gamma_1}(a, \, b) \, = \, \mathfrak{P}_{\Gamma_2}(a, \, b)$ for all but finitely many pairs 
$(a, \, b) \, \in \, \R \, \times \, [0, \, 2\pi]$.
Then $PL_{\Gamma_1}(l) \,= \,  PL_{\Gamma_2}(l)$, \end{theorem}
\medskip

Let $A_{\Gamma_1}^{(0,0)}$ be the multiset on the left hand side of the above equation and 
let $A_{\Gamma_2}^{(0,0)}$ be the multiset on the right hand side.  
Consider now the smallest point $y$ in $A_{\Gamma_1}^{(0,0)}$ such that $y \, > \, 0$.  
This point is of the form $\dfrac{2n_1\pi}{a(p)}$ for some $n_1 \, \in \, \Z$ and $p \, \in \, S_1$.  
Recall that $a(p) \, > \, 0$ for all $p \,\in \,  P_{\Gamma_i}$ ($ i \, = \, 1, \, 2 $).  
If $n_1 \, \leq \, 0$ then $y \, \leq \, 0$.  
Also $\dfrac{2n_1\pi}{a(p)} \, > \, \dfrac{2\pi}{a(p)}$ for $n_1 \, > \, 1$.
Therefore we conclude that $n_1 \, = \, 1$ and $y \, = \, \dfrac{2\pi}{a(p)}$ for some $p \, \in \, S_1$.  
Similarly, we can argue that the smallest point in $A_{\Gamma_2}^{(0,0)}$ is of the form $\dfrac{2\pi}{a(q)}$ 
for some $q \, \in \, S_2$.
The quality of multisets above implies that $\dfrac{2\pi}{a(p)} \, = \, \dfrac{2\pi}{a(q)}$ 
and they occur with same multiplicities.
By removing all the points of the form $\dfrac{2n_1\pi}{a(p)} \; \; \left( \text{resp.} \; \; \dfrac{2n_2\pi}{a(q)}\right)$ from 
$A_{\Gamma_1}^{(0,0)}$ (resp. $A_{\Gamma_2}^{(0,0)}$), we can repeat the argument above to conclude the below equality of sets with multiplicity:

\begin{equation} 
\bigg\{ a(p) \, \mid \, p \, \in \, S_1 \bigg\} \, = 
\bigg\{ a(q) \, \mid \, q \, \in \, S_2  \bigg\}.
\end{equation}\\

Therefore we get the equality of the primitive length spectrum with multiplicity:
\begin{equation}\label{lengtheq}\bigg\{ a(p) \, \mid \, p \, \in \, P_{\Gamma_1} \bigg\} \, = 
\bigg\{ a(q) \, \mid \, q \, \in \, P_{\Gamma_2}  \bigg\}.
\end{equation}\\

Now we look again the case where $\tau$ is the 3-dimensional representation of $\SO(3,1)^{\circ}$ and hence $m \, = \, 1$.
Using the equality of sets in equation \ref{lengtheq}, we can remove the elements of the form $\dfrac{2n_1\pi}{a(p)}$ 
($n_1 \, \in \, \Z$ and $p \,\in \, S_1$) from the set on the left hand side of equation \ref{eqzeroline} and of the form
$\dfrac{2n_2\pi}{a(q)}$ ($n_2 \, \in \, \Z$ and $q \,\in \, S_2$) from the set on the right hand side 
to get the below equality of multisets: \\

\begin{multline}
\label{res2eq}
\bigg\{ \dfrac{-b(p)k_1 - 2n_1\pi}{a(p)} \, 
\mid \, p \, \in \, S_1, \, k_1 \, = \pm1,  \, n_1 \, \in \, \Z \bigg\} \, = \\[2ex]
\bigg\{ \dfrac{-b(q)k_2 - 2n_2\pi}{a(q)} \, 
\mid \, q \, \in \, S_2, \, k_2 \, = \, \pm1,   \, n_2 \, \in \, \Z \bigg\}\
\end{multline} \\

Consider the point in the set on the left hand side which is positive and closest to $0$.  
It is of the form $$s_2 \, = \, \dfrac{-b(p)k \, - \, 2n\pi}{a(p)}$$ 
for some $n \, \in \, \Z$, $p \, \in \, S_1$, and $k \, = \pm 1$.  
We know that $b(p^{-1}) \, = \, 2\pi \, - \, b(p)$ and $a(p^{-1}) \, = \, a(p)$.  Therefore 
\begin{align*}
\dfrac{-kb(p^{-1} \, - \, 2n\pi)}{a(p^{-1})} \, &= \, \dfrac{-k(2\pi \, - \, b(p)) \, - \, 2n\pi}{a(p)} \\
&= \dfrac{b(p)k \, - \, 2(n \, + \, k)\pi}{a(p)}.  
\end{align*}
Replacing $p$ with $p^{-1}$ if needed, we can assume that $k \, = \, -1$.  
Therefore $$s_2 \, = \, \dfrac{b(p) \, - \, 2n\pi}{a(p)}.\\$$

If $n \, > \, 0$, $b(p) \, - \, 2n\pi \, < \, 0$.  But since $s_2 \, > \, 0$, this is not possible and hence 
$n \, \leq \, 0$.
For ease of notation, we put $m \, = \, -n$, which gives 
$$s_2 \, = \, \dfrac{b(p) \, + \, 2m\pi}{a(p)}.$$
But $$ \dfrac{b(p) \, + \, 2m\pi}{a(p)} \, \geq \, \dfrac{b(p)}{a(p)}.$$\medskip

Hence the point in the set on the left hand side which is positive and closest to $0$ is of the form 
$s_2 \, = \, \dfrac{b(p)}{a(p)}$ for some $p \, \in \, S_1$.  
Using similar arguments, one can conclude that 
the point in the set on the left hand side which is positive and closest to $0$ is of the form 
$\dfrac{b(q)}{a(q)}$ for some $q \, \in \, S_2$.  
From the equality of the sets \ref{res2eq}, we get that $\dfrac{b(p)}{a(p)} \, = \, \dfrac{b(q)}{a(q)}$.
Using arguments similar to the ones used in the proof of Theorem \ref{res1}, we get the following theorem: 
\medskip

\begin{theorem}
\label{res2}
Let $G \, = \, \SO(3,1)^{\circ}$ and $\Gamma_1$ and $\Gamma_2$ be two uniform lattices in $G$ such that 
$\mathfrak{P}_{\Gamma_1}(a, \, b) \, = \, \mathfrak{P}_{\Gamma_2}(a, \, b)$ for all but finitely many pairs 
$(a, \, b) \, \in \, \R \, \times \, [0, \, 2\pi]$.
Then $\mathfrak{M}_{\Gamma_1}(c) \,= \,  \mathfrak{M}_{\Gamma_2}(c)$ for all $c \, \in \, \R$.
\end{theorem}\medskip

\begin{remark}
In \cite{BR}, the first author of this paper and Rajan C. S. proved a multiplicity one property for length spectra of even dimensional compact hyperbolic spaces. The methods used there are insufficient to prove an analogous result for odd dimensional spaces due to a different shape of functional equations for Zeta functions. The approach used in this paper is different and we hope to generalise it to higher dimensional hyperbolic spaces.
\end{remark}

\end{document}